\documentclass{amsart}
\usepackage{amscd,amssymb}
\theoremstyle{plain}
\newtheorem{prop}[subsection]{Proposition}
\newtheorem{thm}[subsection]{Theorem}
\newtheorem{lem}[subsection]{Lemma}
\newtheorem{cor}[subsection]{Corollary}

\theoremstyle{remark}
\newtheorem{rem}[subsection]{Remark}
\newtheorem*{ack}{Acknowledgment}

\theoremstyle{definition}
\newtheorem{exm}[subsection]{Example}

\newcommand{\A}{{\mathcal A}}
\newcommand{\B}{{\mathcal B}}
\newcommand{\LL}{{\mathcal L}}
\newcommand{\RR}{{\mathcal R}}

\newcommand{\Z}{{\mathbb Z}}
\newcommand{\C}{{\mathbb C}}
\newcommand{\CP}{{\mathbb{CP}}}
\newcommand{\N}{{\mathbb N}}
\newcommand{\ZN}{{\mathbb Z}_N}
\newcommand{\Q}{{\mathbb Q}}
\newcommand{\T}{{({\mathbb C}^*)^n}}

\newcommand{\UU}{{\mathbf U}}
\newcommand{\VV}{{\mathbf V}}
\newcommand{\WW}{{\mathbf W}}
\newcommand{\bone}{{\mathbf 1}}

\newcommand{\D}{{\Delta}}
\newcommand{\la}{{\lambda }}
\newcommand{\bul}{{\bullet }}

\renewcommand{\a}{{\alpha }}

\renewcommand{\c}{{\gamma }}
\renewcommand{\L}{{\Lambda }}
\renewcommand{\ll}{{\ell }}

\DeclareMathOperator{\rank}{rank}
\DeclareMathOperator{\codim}{codim}
\DeclareMathOperator{\ii}{i}
\DeclareMathOperator{\id}{id}
\DeclareMathOperator{\Mat}{Mat}
\DeclareMathOperator{\TT}{T}

\begin{document}

\title[Arrangements and local sytems]
{Arrangements and local systems}
\author[D.~Cohen]{Daniel C.~Cohen$^\dag$}
\address{Department of Mathematics, Louisiana State University, 
Baton Rouge, LA 70803}
\email{cohen@math.lsu.edu}
\urladdr{http://math.lsu.edu/\~{}cohen}
\thanks{{${\,}^\dag$}Partially supported by 
a grant from the Louisiana State University Council on Research
and by Louisiana Board of Regents grants
LEQSF(1996-99)-RD-A-04 and LEQSF(1999-2002)-RD-A-01}

\author[P.~Orlik]{Peter Orlik$^\ddag$}
\address{Department of Mathematics, University of Wisconsin, 
Madison, WI 53706}
\email{orlik@math.wisc.edu}
\thanks{{${\,}^\ddag$}Partially supported by the NSF}



\begin{abstract}
We use stratified Morse theory to construct a complex to compute the 
cohomology of the complement of a hyperplane arrangement with 
coefficients in a complex rank one local system.  The linearization of 
this complex is shown to be the Orlik-Solomon algebra with the 
connection operator.  Using this result, we establish the relationship 
between the cohomology support loci of the complement and the 
resonance varieties of the Orlik-Solomon algebra for any arrangement, 
and show that the latter are unions of subspace arrangements in 
general, resolving a conjecture of Falk.  We also obtain lower bounds 
for the local system Betti numbers in terms of those of the 
Orlik-Solomon algebra, recovering a result of Libgober and Yuzvinsky.  
For certain local systems, our results provide new combinatorial upper 
bounds on the local system Betti numbers.  These upper bounds enable 
us to prove that in non-resonant systems the cohomology is 
concentrated in the top dimension, without using resolution of 
singularities.
\end{abstract}

\maketitle

\section{Introduction}
\label{sec:intro}

Let $\A=\{H_1,\dots,H_n\}$ be a hyperplane arrangement in $\C^\ll$, 
with complement $M=M(\A)=\C^\ll\setminus\bigcup_{j=1}^n H_j$.  We 
assume that $\A$ contains $\ll$ linearly independent hyperplanes.  Let 
$\la=(\la_1,\dots,\la_n)\in\C^n$ be a collection of weights.  
Associated to $\la$, we have a rank one representation 
$\rho:\pi_1(M)\to\C^*$ given by $\c_j\mapsto t_j=\exp(2\pi\ii\la_j)$ 
for any meridian loop $\c_j$ about the hyperplane $H_j$ of $\A$, and a 
corresponding rank one local system $\LL$ on $M$.  The need to 
calculate the local system cohomology $H^*(M;\LL)$ arises in several 
problems: the Aomoto-Gelfand theory of multivariable hypergeometric 
integrals \cite{AK,Gel1}; representation theory of Lie algebras and 
quantum groups and solutions of the Knizhnik-Zamolodchikov 
differential equation in conformal field theory \cite{Va}; determining 
the cohomology groups of the Milnor fiber of the non-isolated 
hypersurface singularity at the origin obtained by coning the 
arrangement \cite{CS1}.

We call a system of weights $\la\in\C^n$ {\em non-resonant} if the 
Betti numbers of $M$ with coefficients in the associated local system 
$\LL$ are minimal.  The set of non-resonant weights is open and dense 
in $\C^n$.

\begin{thm} 
\label{nonres}
If $\la$ is non-resonant, then 
\[ 
H^q(M;\LL)=0 \mbox{ for } q\neq \ll, \mbox{ and }
\dim H^{\ll}(M;\LL)=|e(M)|.
\]
where $e(M)$ is the Euler characteristic of the complement.
\end{thm}

We use the notation and results of \cite{OT1}.  Let $A=A(\A)$ be the 
Orlik-Solomon algebra of $\A$ generated by the 1-dimensional classes 
$a_j$, $1\leq j\leq n$.  It is the quotient of the exterior algebra 
generated by these classes by a homogeneous ideal, hence a finite 
dimensional graded $\C$-algebra.  There is an isomorphism of graded 
algebras $H^*(M;\C) \simeq A(\A)$.  In particular, $\dim 
A^q(\A)=b_q(\A)$ where $b_q(\A)=\dim H^q(M;\C)$ denotes the $q$-th 
Betti number of $M$ with trivial local coefficients $\C$.
 
Theorem \ref{nonres} is a consequence of two results which are quite 
different in nature.  One is the work of Esnault-Schechtman-Viehweg 
\cite{ESV}, refined by Schechtman-Terao-Varchenko \cite{STV}, and 
obtained by using Deligne's work \cite{De}.  It involves resolution of 
singularities and techniques of algebraic geometry.  Here the de Rham 
complex with a twisted differential is used to compute $H^*(M;\LL)$, 
and the main aim is to reduce the infinite dimensional cochain groups 
to a complex with finite dimensional cochain groups.  In the passage 
from consideration of global rational differential forms in $\C^\ll$ 
with poles of arbitrary order on the hyperplanes of $\A$ to poles of 
order one (logarithmic), certain geometric conditions on $\la$ arise.  
The final result is that for suitable $\la$ there is an isomorphism
\[
H^q(M;\LL)\simeq H^q(A^{\bul},a_\la\wedge),
\]
where the Orlik-Solomon algebra with differential $A^q\to A^{q+1}$ 
given by multiplication by $a_\la=\sum_{j=1}^{n}\la_j\,a_j$ is 
identified with a subcomplex of the twisted de Rham complex.  The 
second result, due to Yuzvinsky \cite{Yuz}, is purely combinatorial in 
nature.  It shows that if $\la$ satisfies certain combinatorial 
conditions, then
\[ 
H^q(A^{\bul},a_\la\wedge)=0 \mbox{ for } q\neq \ll, \mbox{ and }
\dim H^{\ll}(M;\LL)=|e(M)|.
\]

The geometric conditions on $\la$ noted above arise from conditions on 
the monodromy of the local system about certain intersections of 
hyperplanes of $\A$.  These conditions derive from Aomoto's 1973 work 
\cite{Ao}, where slightly stronger monodromy conditions were 
considered.  Indeed, the results of \cite{ESV} mentioned previously 
resolved a conjecture from this paper.  Kohno \cite{Ko} used similar 
(strong) monodromy conditions to prove a vanishing theorem analogous 
to Theorem~\ref{nonres} for local system cohomology and locally finite 
homology, and exhibited a basis for the non-zero homology group in a 
special case.  Hattori \cite{Ha} showed that for general position 
arrangements only the trivial local system is resonant.  More 
recently, Falk-Terao \cite{FT} constructed a basis for the non-zero 
local system cohomology group in the non-resonant case for any 
arrangement.  Thus, much is known in the case of non-resonant weights.  
See \cite{OT2} for a detailed discussion.
 
Substantially less is known about resonant weights, which are of 
particular interest in several of the applications noted above.  This 
interest has been a motivating factor in much recent work, including 
\cite{C3,CS1,CS4,Fa,Li,LY}.  In this paper, the local system 
cohomology is studied using a complex to compute $H^{*}(M;\LL)$ which 
arises from stratified Morse theory.  The terms of this complex are 
finite dimensional but the boundary maps are not easily computed.  In 
fact, we construct a universal complex, 
$(K^{\bul}_{\L}(\A),\D^{\bul}(x))$, where $x=(x_1,\ldots,x_n)$ are 
non-zero complex variables, with the property that the specialization 
$x_j \mapsto t_j$ calculates $H^*(M;\LL)$.  There is a similar 
universal complex, called the Aomoto complex 
$(A^{\bul}_{R}(\A),a_{y}\wedge)$, where $y=(y_1,\ldots,y_n)$ are 
variables, with the property that the specialization $y_j \mapsto 
\la_j$ calculates $H^*(A^{\bul},a_\la\wedge)$.  Our first result is:

\begin{thm} 
\label{approx}
For any arrangement $\A$, the Aomoto complex 
$(A^{\bul}_{R}(\A),a_{y}\wedge)$ is chain equivalent to the 
linearization of the universal complex 
$(K^{\bul}_{\L}(\A),\D^{\bul}(x))$.
\end{thm}

Since this theorem applies to all weight systems, we obtain 
corollaries for resonant weights.  One is a proof of Falk's conjecture 
that the resonance varieties of the Orlik-Solomon algebra are unions 
of subspace arrangements, see Corollary \ref{cor:falkconj}.  Another 
is the known lower bound of Proposition \ref{prop:lb}:
\[
\dim H^q(A^{\bul},a_\la\wedge)\leq \dim H^q(M;\LL).
\]

In the other direction, the upper bounds 
\[
\dim H^q(M;\LL)\leq \dim H^q(M;\C)
\]
were conjectured in \cite{AK}, and established in \cite{C2}.  In most 
applications where resonant weights occur, the $\la_j$ are rational 
numbers.  We call the associated local system {\em rational}.  Our 
next aim is to obtain better upper bounds for $\dim H^q(M;\LL)$ for 
rational local systems.  Let $\la_j= k_j/N$, where $k_j\in\Z$ and 
$N\in\N$.  We define a complex, $(A^{\bul}_N,\bar{a}_k\wedge)$, whose 
cochain groups are the graded parts of the Orlik-Solomon algebra with 
coefficients in the ring $\ZN=\Z/N\Z$, and whose differential operator 
is obtained from $a_k=\sum_{j=1}^{n}k_j\,a_j$ by reduction mod $N$, 
and prove:

\begin{thm}
Let $\la=k/N$ be a system of rational weights, and let $\LL$ be the 
associated rational local system on the complement $M$ of $\A$.  
Then, for each $q$,
\begin{equation*} 
\dim_\C H^q(M;\LL) \le \rank_{\ZN} H^q(A^{\bul}_N,\bar{a}_k\wedge).
\end{equation*}
\end{thm} 

This result about rational, and often resonant, weights is used to 
give a proof of Theorem \ref{nonres} (about non-resonant weights) 
which does not rely on resolution of singularities.  Examples show 
that the inequalities in both the upper and lower bounds may be 
strict.

\section{Linear Approximation}

Let $\A=\{H_1,\dots,H_n\}$ be a hyperplane arrangement in $\C^\ll$, 
and let $\LL$ be a complex rank one local system on the complement $M$ 
of $\A$.  In this section, we show that the Orlik-Solomon algebra 
complex $(A^\bul(\A),a_\la\wedge)$ is, in a natural sense, a linear 
approximation of a Morse theoretic complex, $(K^\bul(\A),\D^\bul)$, 
the cohomology of which is isomorphic to that of $M$ with coefficients 
in $\LL$.

\subsection{Preliminaries}
Let $\B$ denote the Boolean arrangement in $\C^n$, with coordinates 
$z_1,\dots,z_n$.  A defining polynomial for $\B$ is given by 
$Q(\B)=z_1\cdots z_n$, and the complement of $\B$ is the complex 
$n$-torus, $T=\T$.

\begin{prop} \label{prop:slice}
Let $\A$ be an arrangement of $n$ hyperplanes in $\C^\ll$.  Then the 
complement $M$ of $\A$ may be realized as a linear section of the 
complex $n$-torus $T$.
\end{prop}
\begin{proof}
Without loss of generality, assume that $\A$ is an essential 
arrangement in $\C^\ll$, so that $\A$ contains $\ll$ independent 
hyperplanes.  Choose coordinates $z_1,\dots,z_\ll$ for $\C^\ll$.  
Ordering the hyperplanes of $\A$ appropriately, a defining polynomial 
for $\A$ is then given by $Q(\A)=z_1\cdots 
z_\ll\cdot\a_{\ll+1}\cdots\a_n$, where $\a_j=\a_j(z_1,\dots,z_\ll)$ is 
a linear polynomial.

Let $S$ denote the $\ll$-dimensional affine subspace of $\C^n$ defined 
by the equations $\{z_j=\a_j(z_1,\dots,z_\ll) \mid \ll+1\le j\le n\}$.  
Identifying $\C^\ll$ with $S$ in the obvious manner, it is readily 
checked that $M = S\cap T$.
\end{proof}

Let $\la=(\la_1,\dots,\la_n)\in\C^n$ be a collection of weights.  
Associated to $\la$, we have a rank one representation 
$\rho:\pi_1(M)\to\C^*$ given by $\c_j\mapsto t_j=\exp(2\pi\ii\la_j)$ 
for any meridian loop $\c_j$ about the hyperplane $H_j$ of $\A$, and a 
corresponding rank one local system $\LL$ on $M$.  Note that the 
representation $\rho$ factors through the first homology $H_1(M;\Z)$, 
which is generated by the classes $[\c_j]$.  Since these classes also 
generate $H_1(T;\Z)=\pi_1(T)$, the local system $\LL$ extends to a 
local system on $T$, which we continue to denote by $\LL$.

\subsection{The Morse theoretic complex}
For any complex local system $\LL$ on the complement of an arrangement 
$\A$, in \cite{C1} we used stratified Morse theory to construct a 
complex $(K^\bul(\A),\D^\bul)$, the cohomology of which is naturally 
isomorphic to $H^*(M;\LL)$, the cohomology of $M$ with coefficients in 
$\LL$.  We give a brief description of this complex (for a rank one 
local system) and record some relevant results from \cite{C1}.

An {\em edge} of $\A$ is a nonempty intersection of hyperplanes.  The 
arrangement $\A$ determines a Whitney stratification of $S=\C^\ll$, 
with a stratum, $S_X=X \setminus \bigcup_{Y\subsetneq X} Y$, of 
codimension $p$ associated to each codimension $p$ edge $X$ of $\A$.  
(The stratum $S_X$ is given by $M(\A^X)$, see \cite{OT1}.)  Let 
$\emptyset = F_{-1} \subset F_0 \subset F_1 \subset F_2 \subset\dots 
\subset F_\ll = S$ be a flag in $S$ which is transverse to the 
stratification determined by $\A$, so that $\dim F_q\cap S_X = 
q-\codim S_X$ for each stratum, where a negative dimension indicates 
that $F_q\cap S_X=\emptyset$.  Such a flag may be constructed using a 
Morse function on $S$ that is weakly self-indexing with respect to the 
above stratification, see \cite[Section~1]{C1}.

The sets $M_q=F_q\cap M$ form a filtration, $\emptyset = M_{-1} 
\subset M_0 \subset M_1 \subset \dots \subset M_\ll = M$.  By 
construction, for each $q$, the closure of ${M}_q$ intersects all 
strata of codimension at most $q$, and intersects no stratum of 
codimension greater than $q$.  This filtration is well-suited for the 
study of local system cohomology in the sense of the following.

\begin{prop} \label{prop:vanish}
For each $q$, $0\le q \le \ll$, we have $H^i(M_q,M_{q-1};\LL) = 0$ if 
$i \neq q$, and $\dim H^q(M_q,M_{q-1};\LL) = b_q(\A)$. 
\end{prop}
This Proposition may be proved using stratified Morse theory 
\cite{GM}.  For details, the reader is referred to \cite[Sections~2, 
3, and 5]{C1}.

For each $q$, let $K^q(\A)=H^q(M_q,M_{q-1};\LL)$ and denote by 
$\D^q$ the boundary homomorphism $H^{q}(M_{q},M_{q-1};\LL) \to 
H^{q+1}(M_{q+1},M_{q};\LL)$ of the triple $(M_{q+1},M_q,M_{q-1})$.  
It 
is readily checked that the composition $\D^{q+1}\circ\D^q=0$.  Thus 
the system of complex vector spaces and linear maps 
$(K^\bul(\A),\D^\bul)$ is a complex.  The following is a special case 
of \cite[Theorem~2.4]{C1}.

\begin{thm} \label{thm:Kdot} 
The cohomology of the complex $(K^\bul(\A),\D^\bul)$ is naturally 
isomorphic to $H^*(M;\LL)$, the cohomology of $M$ with coefficients in 
the local system $\LL$.
\end{thm}

\begin{rem} As shown in \cite{C2}, the inequalities $\dim 
H^q(M;\LL)\le \dim H^q(M;\C)$ noted in the Introduction  
follow immediately from this result.  See also \cite{Ma}.
\end{rem}

\begin{exm} \label{exm:boolean}
We discuss first the important special case of the torus.  Recall that 
$\B$ denotes the Boolean arrangement, consisting of the coordinate 
hyperplanes in $\C^n$, and that a rank one local system $\LL$ on the 
complement of any arrangement of $n$ hyperplanes extends naturally to 
a local system on the torus $T=\T$, the complement of $\B$.  For this 
arrangement (and more generally any general position arrangement), the 
Morse theoretic complex constructed above admits a complete 
description, see \cite[Section~7]{C1} for details.  Denote this 
complex by $(K^\bul(\B),D^\bul)$.

By Proposition~\ref{prop:vanish} above, the terms of this complex are 
vector spaces of dimension $\dim K^q(\B) = \dim H^q(T;\C) = b_q(\B) 
= \binom{n}{q}$.  Thus, as a graded group, the complex $K^\bul(\B)$ 
may be identified with the cohomology of $T$, which in turn is 
isomorphic to the exterior algebra $E=\bigwedge E^1$, where 
$E^1=\oplus_{j=1}^n \C{e_j}$.

The local system $\LL$ is induced by the representation 
$\rho:\pi_1(T)\to \C^*$ defined by $\c_j\mapsto t_j$, where $\c_j$ is 
a meridian loop about the hyperplane $H_j$ and 
$t_j=\exp(2\pi\ii\la_j)$.  Under the identification $K^\bul(\B)=E$ 
above, the boundary map $D^{q}:K^q(\B) \to K^{q+1}(\B)$ is given 
by
\begin{equation*} \label{eq:etwedge}
D^{q}(e_J)=e_t \wedge e_J=\bigl((t_1-1)e_1+\dots+(t_n-1)e_n\bigr)
\wedge e_J,
\end{equation*}
where $e_J=e_{j_1}\wedge\dots\wedge e_{j_q}$ if $J=\{j_1,\dots,j_q\}$.
\end{exm}

\subsection{Approximation}  
We now return to an arbitrary arrangement $\A$ of $n$ hyperplanes and 
show that the Morse theoretic complex described above is approximated 
by the Orlik-Solomon algebra complex.  Let $\la=(\la_1,\dots,\la_n)$ 
be a weight vector in $\C^n$, $t=(t_1,\dots,t_n)$ the associated point 
in the complex torus $\T$, and $\LL$ the corresponding local system on 
$M(\A)$.  It is important to note that if $m=(m_1,\ldots,m_n)$ is a 
tuple of integers, then $\la +m=(\la_1+m_1,\dots,\la_n+m_n)$ 
determines the same $t$ and $\LL$.

As evidenced by Proposition~\ref{prop:vanish}, the dimensions of the 
terms of the complex $(K^\bul(\A),\D^\bul)$ are independent of $t$ 
(resp.,~$\la$).  To indicate the dependence of the complex on $t$, we 
write $\D^\bul=\D^\bul(t)$, $\LL=\LL_t$ and view these boundary maps 
as functions of $t$.  These observations may be formalized as follows.  
Let $\L=\C[x_1^{\pm 1},\dots,x_n^{\pm 1}]$ denote the ring of complex 
Laurent polynomials in $n$ commuting variables, and for each $q$, let 
$K^q_\L(\A)= \L \otimes_{\C} K^q(\A)$.

\begin{thm} \label{thm:univcx}
Given an arrangement $\A$ of $n$ hyperplanes with complement $M$, 
there exists a universal complex $(K^\bul_\L(\A),\D(x))$ with the 
following properties:
\begin{enumerate}
\item The terms are free $\L$-modules, whose ranks are given by the 
Betti numbers of $M$, $K^q_\L(\A) \simeq \L^{b_q(\A)}$.

\item The boundary maps, $\D^q(x): K^q_\L(\A) \to K^{q+1}_\L(\A)$ are 
$\L$-linear.

\item For each $t\in\T$, the specialization $x \mapsto t$ yields the 
complex $(K^\bul(\A),\D^\bul(t))$, the cohomology of which is 
isomorphic to $H^*(M;\LL_t)$, the cohomology of $M$ with coefficients 
in the local system associated to $t$.
\end{enumerate}
\end{thm}

\begin{rem} \label{rem:boolean}
For the Boolean arrangement $\B$ of $n$ hyperplanes, the universal 
complex $(K^\bul_\L(\B),D^\bul(x))$ is dual to the standard free 
$\L\simeq\Z\Z^n$-resolution of the integers.  This follows from the 
description of the complex $(K^\bul(\B),D^\bul(t))$ given in 
Example~\ref{exm:boolean}.
\end{rem}

\begin{prop} \label{prop:univchain}
Let $\A$ be an arrangement of $n$ hyperplanes in $\C^{\ll}$.  Then 
there exists a chain map $\Psi^{\bul}(x):K^\bul_\L(\B) \to 
K^\bul_\L(\A)$ from the universal complex of the Boolean arrangement 
$\B$ to the universal complex of $\A$.
\end{prop}
\begin{proof}
Recall from Proposition~\ref{prop:slice} that we realize the 
complement of $\A$ as an $\ll$-dimensional linear section, $M=S \cap 
T$, of the complement of $\B$, and that any rank one local system 
$\LL$ on $M$ extends to a local system on $T$.

For such a local system, the complexes $K^{\bul}(\A)$ and $K^\bul(\B)$ 
are constructed using flags $\emptyset= F_{-1} \subset F_0 \subset F_1 
\subset \dots \subset F_\ll=S$ and $\emptyset= F'_{-1} \subset F'_0 
\subset F'_1 \subset \dots \subset F'_n=\C^n$ respectively.  Since the 
affine subspace $S$ is transverse to the hyperplanes of $\B$, we may 
assume that $F'_q=F_q$ for $q\le 1$.  Thus the inclusion $M\cap 
F_q\subseteq T\cap F_q$ induces canonical isomorphisms $\psi^q:K^q(\B) 
\xrightarrow{\sim} K^q(\A)$ for $q\le 1$, and for every $t\in\T$ we 
have $\D^0(t)=D^0(t)$.  This yields isomorphisms $\Psi^{q}=\id \otimes 
\psi^{q}:K^{q}_{\L}(\B) \xrightarrow{\sim} K^{q}_{\L}(\A)$ for $q\le 
1$, and we have $\D^{0}(x)=D^{0}(x)$.

Now, as noted in Remark~\ref{rem:boolean}, the complex 
$(K^{\bul}_{\L}(\B),D^{\bul}(x))$ is dual to the standard 
$\L$-resolution of $\Z$, which is of course acyclic.  Thus it follows 
from the acyclic models theorem that there is a chain map 
$\Psi^{\bul}(x): K^{\bul}_{\L}(\B) \to K^{\bul}_{\L}(\A)$ covering 
$\Psi^{q}$ ($q\le 1$).
\end{proof}

As is the case for the local system cohomology of the complement, 
there is a universal complex for the cohomology, 
$H^{*}(A^{\bul}(\A),a_{\la}\wedge)$, of the Orlik-Solomon algebra.  
Let $R=\C[y_{1},\dots,y_{n}]$ be the polynomial ring.  The {\em Aomoto 
complex} $(A^{\bul}_{R}(\A),a_{y}\wedge)$ has terms 
$A^{q}_{R}(\A)=R\otimes_{\C} A^{q}(\A)$, and boundary maps given by 
$p(y)\otimes\eta \mapsto \sum y_{j}p(y) \otimes a_{j} \wedge \eta$.  
For $\la \in \C^{n}$, the specialization $y\mapsto \la$ of the Aomoto 
complex $(A^{\bul}_{R}(\A),a_{y}\wedge)$ yields the Orlik-Solomon 
algebra complex $(A^{\bul}(\A),a_{\la}\wedge)$.

Fix a basis for the Orlik-Solomon algebra of $\A$.  Evidently, this 
yields a basis for each term $A^{q}_{R}(\A)$ of the Aomoto complex.  
Let $\mu^{q}(y)$ denote the matrix of $a_{y}\wedge:A^{q}_{R}(\A) \to 
A^{q+1}_{R}(\A)$ with respect to this basis.

\begin{lem} \label{lem:linear}
For each $q$, the entries of $\mu^{q}(y)$ are integral linear forms 
in $y_1,\dots,y_n$.
\end{lem}
\begin{proof}
First note that this holds in the case where $\A=\B$ is the Boolean 
arrangement.  In this case, the Orlik-Solomon algebra is the exterior 
algebra $E$, and the boundary maps of the Aomoto complex are given by 
$e_J \mapsto e_y\wedge e_J= \sum_{j=1}^n y_j\otimes e_j\wedge e_J$ on 
generators.

For an arbitrary arrangement $\A$, the Aomoto complex 
$(A^{\bul}_{R}(\A),a_{y}\wedge)$ may be realized as the quotient of 
the Aomoto complex $(A^{\bul}_{R}(\B),e_{y}\wedge)$ of the Boolean 
arrangement by the subcomplex $(I^{\bul}_{R}(\A),e_{y}\wedge)$, where 
$I^\bul_R(\A)$ denotes the tensor product of $R$ with the 
Orlik-Solomon ideal $I(\A)$.  Since the ideal $I(\A)$ is defined by 
integral linear combinations of the generators of the exterior 
algebra, the result follows.
\end{proof}

The main result of this section is the following.

\begin{thm} \label{thm:approx}
For any arrangement $\A$, the Aomoto complex 
$(A^{\bul}_{R}(\A),a_{y}\wedge)$ is chain equivalent to the 
linearization of the universal complex
$(K^{\bul}_{\L}(\A),\D^{\bul}(x))$
\end{thm}
\begin{proof}
For each $q$, fix a basis for $K^{q}_{\L}(\A)$.  Since $\rank_{\L} 
K^{q}_{\L}(\A) = \dim_{\C} A^{q}(\A)=b_{q}(\A)$, the basis for 
$K^{\bul}_{\L}(\A)$ may be chosen in one-to-one correspondence with 
that of $A(\A)$.  We shall not distinguish between the boundary map, 
$\D^{q}(x):K^{q}_{\L}(\A) \to K^{q+1}_{\L}(\A)$, of the universal 
complex and its matrix with respect to the chosen basis.

The entries of $\D^{q}(x)$ are elements of the Laurent polynomial ring 
$\L$, the coordinate ring of the complex algebraic $n$-torus.  Via the 
specialization $x \mapsto t \in \T$, we shall view them as holomorphic 
functions $\T\to\C$.  Similarly, for each $q$, we view $\D^{q}$ as a 
holomorphic map $\D^{q}:\T\to\Mat(\C)$ by $t\mapsto \D^{q}(t)$.

Let $\bone=(1,\dots,1)$ denote the identity element of $\T$.  The 
holomorphic tangent space of the complex $n$-torus at $\bone$ is 
$\TT_{\bone}\T = \C^{n}$.  Identify the coordinates 
$y=(y_{1},\dots,y_{n})$ of this tangent space with the variables 
appearing in the Aomoto complex.  The exponential map $\TT_{\bone}\T 
\to \T$ is induced by $\exp:\C\to\C^{*}$, $y_{j}\mapsto e^{y_{j}} = 
x_{j}$.

The specialization $x \mapsto \bone$ corresponding to the trivial 
local system yields a complex $(K^{\bul}(\A),\D^{\bul}(\bone))$ with 
trivial boundary maps, $\D^{q}(\bone)=0$ for each $q$.  This follows 
from Proposition~\ref{prop:vanish}.  See \cite[III.3]{GM} for a 
detailed discussion of this phenomenon in homology.

If $f \in \L$ is an entry of $\D^{q}(x)$, then the map $f:\T\to\C$ 
satisfies $f(\bone)=0$.  For such a Laurent polynomial (resp.,~map), 
the derivative at the identity, $f_{*}:\TT_{\bone}\T = \C^{n} \to \C= 
\TT_{0}\C^{*}$, is given by
\[
f_*(y)={\frac{d}{ds}} f(e^{s y_1},\dots,e^{s y_n})\bigr\rvert_{s=0}.
\]  
Similarly, for $F\in \Mat(\L)$ satisfying $F(\bone)=0$, we have 
$F_{*}:\TT_{\bone}\T \to \TT_{0}\Mat(\C)=\Mat(\C)$.  For two Laurent 
polynomials $f$ and $g$, the product rule yields $(fg)_*(y)=f_*(y) 
g(\bone)+ f(\bone) g_*(y)$.  More generally, for $F \in \Mat_{p\times 
q}(\L)$ and $G \in \Mat_{q\times r}(\L)$, using matrix multiplication 
and the differentiation rules we have
\begin{equation}  
\label{eq:leibniz}
(FG)_*(y)=F_*(y)\cdot G(\bone)+F(\bone)\cdot G_*(y).
\end{equation}

Now recall the chain map $\Psi^{\bul}(x):K^\bul_\L(B) \to 
K^\bul_\L(\A)$ from the universal complex of the Boolean arrangement 
$\B$ to that of $\A$ described in Proposition~\ref{prop:univchain}.  
As is the case for the boundary maps $D^{q}(x)$ and $\D^{q}(x)$ of the 
two complexes, we view $\Psi^{q}(x)$ as a holomorphic map 
$\T\to\Mat(\C)$, and we do not distinguish between this map and its 
matrix.  Since $\Psi^{\bul}(x)$ is a chain map, we have $D^{q}(x) 
\cdot \Psi^{q+1}(x) = \Psi^{q}(x) \cdot \D^{q}(x)$ for each $q$.  We 
differentiate at $\bone$ using the product rule~\eqref{eq:leibniz} 
above to obtain
\[
D^{q}_{*}(y) \cdot \Psi^{q+1}(\bone) + D^{q}(\bone) \cdot 
\Psi^{q+1}_{*}(y) = \Psi^{q}_{*}(y) \cdot \D^{q}(\bone) + 
\Psi^{q}(\bone) \cdot \D^{q}_{*}(y).
\]
Since $D^{q}(\bone)=0$ and $\D^{q}(\bone)=0$ for all $q$, we thus have
\[
D^{q}_{*}(y) \cdot \Psi^{q+1}(\bone) = 
\Psi^{q}(\bone) \cdot \D^{q}_{*}(y).
\]

Recall that $R=\C[y_{1},\dots,y_{n}]$ denotes the polynomial ring. 
For $F\in\Mat(\L)$, view the derivative $F_{*}(y)$ as a linear map 
between free $R$-modules in the obvious way.  Consider the systems of 
free $R$-modules and $R$-linear maps
\[
(K^{\bul}_{R}(\A),\D^{\bul}_{*}(y)) \quad \text{and} \quad
(K^{\bul}_{R}(\B),D^{\bul}_{*}(y)),
\]
where $K^{q}_{R}(\A)=R \otimes_{\C} K^{q}(\A)$ and $K^{q}_{R}(\B)=R 
\otimes_{\C} K^{q}(\B)$.  From the description of the universal 
complex of the Boolean arrangement stemming from 
Example~\ref{exm:boolean} and the fact that $(x_{j})_{*}=y_{j}$ for 
each $j$, it is clear that the system 
$(K^{\bul}_{R}(\B),D^{\bul}_{*}(y))$ is a complex, and coincides with 
the Aomoto complex of the arrangement $\B$.

We assert that the system $(K^{\bul}_{R}(\A),\D^{\bul}_{*}(y))$ is 
also a complex, and is chain equivalent to the Aomoto complex of $\A$.  
For this, consider again the specialization $x\mapsto\bone$ 
corresponding to the trivial local system $\LL=\C$.  In this instance, 
as noted above, the boundary maps of $K^\bul(\A)$ and $K^\bul(\B)$ are 
all trivial.  Thus, for $\LL$ trivial, these complexes simply record 
the cohomology $H^*(M;\C)$ and $H^*(T;\C)$.  As is well known, both 
algebras are generated in dimension one, and the inclusion $M\subset 
T$ induces an epimorphism in cohomology.  Now recall from the proof of 
Proposition~\ref{prop:univchain} that for $q\le 1$, $\Psi^{q}=\id 
\otimes \psi^{q}$ is constant, and that $\psi^{q}:H^{q}(T;\C) \to 
H^{q}(M;\C)$ is induced by inclusion.  So by the naturality of cup 
products and the continuity of $\Psi^{\bul}(x)$, we have 
$\Psi^{q}(\bone)=\id \otimes \psi^{q}$ for all $q$, where 
$\psi^{q}:H^{q}(T;\C)\to H^{q}(M;\C)$ is induced by $\psi^{1}$.

Identifying the cohomology of $T$ with the exterior algebra, 
$H^{*}(T;\C)=E$, and the cohomology of $M$ with the Orlik-Solomon 
algebra, $H^{*}(M;\C)=A(\A)$, we realize the map $\Psi^{\bul}(\bone)$ 
as the map on $R$-modules induced by a choice of projection $\psi:E\to 
A$ from the exterior algebra to the Orlik-Solomon algebra.  
Consequently, the system of $R$-modules and $R$-linear maps 
$(K^{\bul}_{R}(\A),\D^{\bul}_{*}(y))$ is a complex, and is chain 
equivalent to the Aomoto complex of $\A$.
\end{proof}

\section{Some Consequences} \label{sec:apps1}
We discuss some immediate applications of Theorem~\ref{thm:approx}.  

\subsection{Lower Bounds}
We first show that the Orlik-Solomon algebra cohomology provides a 
lower bound for the local system cohomology of the complement, 
recovering a result of Libgober-Yuzvinsky, see \cite[Proposition 4.2, 
Corollary 4.3]{LY}.  Fix an arrangement $\A$ of $n$ hyperplanes, with 
complement $M=M(\A)$ and Orlik-Solomon algebra $A=A(\A)$.  Recall that 
each weight vector $\la\in\C^n$ gives rise to an associated rank one 
local system $\LL=\LL_\la$ on $M$.

\begin{prop}
\label{prop:lb}
For each $\la\in\C^n$ and each $q$, we have 
\begin{equation} \label{eq:LB}
\sup_{m\in\Z^{n}} \dim H^{q}(A^{\bul},a_{\la+m}\wedge) \le 
\dim H^{q}(M;\LL).
\end{equation}
\end{prop}
\begin{proof}
First note that this result clearly holds for $\la=0$.
In general, given $\la=(\la_1,\dots,\la_n)\in\C^n$, the cohomology, 
$H^*(A^{\bul},a_{\la}\wedge)$, of the Orlik-Solomon algebra is 
given by that of the specialization, $y\mapsto\la$, of the Aomoto 
complex $(A^{\bul}_{R}(\A),a_{y}\wedge)$.  Recall that we denote the 
matrix of $a_y\wedge$ by $\mu^q(y)$.  Similarly, the local system 
cohomology, $H^*(M;\LL)$, may be computed from the specialization, 
$x\mapsto t$, of the universal complex 
$(K^{\bul}_{\L}(\A),\D^{\bul}(x))$, where $t=(t_1,\dots,t_n) \in \T$ 
satisfies $t_j=\exp(2\pi\ii\la_j)$, and $\LL$ is induced by the 
representation $\pi_1(M)\to\C^*$, $\c_j\mapsto t_j$ for any meridian 
loop $\c_j$ about $H_j\in\A$.

Now by Theorem~\ref{thm:approx}, the Aomoto complex is chain 
equivalent to $(K^{\bul}_{R}(\A),\D^{\bul}_*(y))$, the linearization 
of the universal complex.  Furthermore, it is well-known that the 
cohomology, $H^*(A^{\bul},a_{\la}\wedge)$, of the Orlik-Solomon 
algebra is invariant under (non-zero) rescaling.  Thus for $\la$ 
sufficiently small, the Inverse Function Theorem implies that
\[
\rank\mu^q(\la)=\rank\D^q_*(\la)=\rank\D^q_*(2\pi\ii\la)=\rank \D^q(t)
\]
for each $q$.  Thus, $\dim H^{q}(A^{\bul},a_{\la}\wedge) = \dim 
H^{q}(M;\LL)$ for $\la$ sufficiently small.  Hence, using invariance 
under rescaling for the cohomology of the Orlik-Solomon algebra and 
upper semicontinuity for the local system cohomology of the 
complement, we have $\dim H^{q}(A^{\bul},a_{\la}\wedge) \le \dim 
H^{q}(M;\LL)$ for arbitrary $\la$.

Now let $m=(m_1,\dots,m_n)\in \Z^n$, and recall that the local system 
associated to $\la+m$ coincides with that associated to $\la$.  Thus, 
by applying the above considerations to the weight vector $\la+m$, we 
see that $\dim H^{q}(A^{\bul},a_{\la+m}\wedge) \le \dim H^{q}(M;\LL)$ 
for each $q$ and all $m\in \Z^n$, which completes the proof.
\end{proof}

\begin{rem} \label{rem:LYsup}
In the case $q=1$, Libgober-Yuzvinsky show that equality holds in
\eqref{eq:LB} for almost all $\la$, see \cite[Theorem 5.3]{LY}.
\end{rem}

There are local systems $\LL$ (resp.,~weight vectors $\la$) for which 
the inequality \eqref{eq:LB} is strict, as illustrated by the 
following examples.  Both involve central arrangements, where all 
local system cohomology groups vanish for non-resonant weights.  Thus 
a weight vector which produces cohomology is a priori resonant.

\begin{exm} \label{exm:lstrict}
Let $\A=\{H_1,\dots ,H_7\}$ be the central arrangement in $\C^{3}$ 
with defining polynomial $Q(\A)=x(x+y+z)(x+y-z)y(x-y-z)(x-y+z)z$.  
(The ordering of the hyperplanes corresponds to that of the factors of 
$Q(\A)$.)  Consider the weight vector 
$\la=\frac{1}{2}(1,0,0,1,1,0,1)\in\C^{7}$, the associated point 
$t=(-1,1,1,-1,-1,1,-1)\in(\C^{*})^{7}$, and the corresponding local 
system $\LL$ on the complement $M$ of $\A$.

An exercise in the Orlik-Solomon algebra reveals that, for all 
$m\in\Z^{7}$, we have $\dim H^{1}(A^{\bul}(\A),a_{\la+m}\wedge)\le 1$.  
On the other hand, it is known that the point $t\in \Sigma^{1}_{2}(M)$ 
is in the second cohomology support locus of (the first cohomology of) 
$M$, see \cite[Example~4.4]{CS4} and Section~\ref{subsec:loci} below.  
It follows that $\dim H^{1}(M(\A),\LL)=2$, and the inequality 
\eqref{eq:LB} is strict for this local system (resp.,~system of 
weights).
\end{exm}

\begin{exm} \label{exm:ceva}
A similar example is provided by the Ceva(3) arrangement $\A$ (the 
monomial arrangement $\A=\A_{3,3,3}$), defined by 
$Q(\A)=(x^3-y^3)(x^3-z^3)(y^3-z^3)$.  For the weight vector 
$\la=\frac{1}{3}(1,1,1,1,1,1,-2,-2,-2)\in\C^9$ and the corresponding 
local system $\LL$ on the complement $M$ of $\A$, it is known that 
$\dim H^{1}(A^{\bul}(\A),a_{\la+m}\wedge)\le 1$ for all $m\in\Z^9$, in 
particular $\dim H^{1}(A^{\bul}(\A),a_{\la}\wedge)=1$, and that $\dim 
H^{1}(M(\A),\LL)=2$.
\end{exm}

\subsection{Cohomology Support Loci and Resonance Varieties}
\label{subsec:loci}
As another application of Theorem~\ref{thm:approx}, we establish the 
relationship between the cohomology support loci of the complement of 
an essential arrangement $\A$ in $\C^{\ll}$ and the resonance 
varieties of its Orlik-Solomon algebra.

Recall that each point $t\in\T$ gives rise to a local system 
$\LL=\LL_{t}$ on the complement $M=M(\A)$.  For non-resonant $t$, the 
cohomology $H^{q}(M,\LL_{t})$ vanishes for $q<\ll$, see Theorem 
\ref{nonres}.  Those $t$ for which $H^{q}(M;\LL_t)$ does not vanish 
comprise the cohomology support loci
\[
\Sigma^{q}_m(M)=\{t\in\T \mid \dim H^{q}(M;\LL_{t})\ge m\}.
\]
These loci are algebraic subvarieties of $\T$, which are invariants of 
the homotopy type of $M$.  See Arapura \cite{Ar} and Libgober 
\cite{Li} for detailed discussions of these varieties in the contexts 
of quasiprojective varieties and plane algebraic curves.

Similarly, each point $\la\in\C^{n}$ gives rise to an element 
$a_{\la}\in A^{1}$ of the Orlik-Solomon algebra $A=A(\A)$.  For 
sufficiently generic $\la$, the cohomology 
$H^{q}(A^{\bul},a_\la\wedge)$ vanishes for $k<\ll$, see \cite{Yuz,Fa}.  
Those $\la$ for which $H^{q}(A^{\bul},a_\la\wedge)$ does not vanish 
comprise the resonance varieties
\[
\RR^{q}_m(A)=\{\la\in\C^{n} \mid \dim H^{q}(A^{\bul},a_\la\wedge)\ge 
m\}.
\]
These subvarieties of $\C^{n}$ are invariants of the Orlik-Solomon 
algebra $A$.  See Falk \cite{Fa} and Libgober-Yuzvinsky \cite{LY} for 
detailed discussions of these varieties.

Recall that $\bone=(1,\dots,1)$ denotes the identity element of $\T$.

\begin{thm} \label{thm:tcone}
Let $\A$ be an arrangement in $\C^\ll$ with complement $M$ and 
Orlik-Solomon algebra $A$.  Then for each $q$ and $m$, the resonance 
variety $\RR_{m}^{q}(A)$ coincides with the tangent cone of the 
cohomology support locus $\Sigma_{m}^{q}(M)$ at the point $\bone$.
\end{thm}
\begin{proof}
For each $t\in\T$, the cohomology of $M$ with coefficients in the 
local system $\LL_{t}$ is isomorphic to that of the Morse theoretic 
complex $(K^\bul(\A),\D^\bul(t))$, the specialization at $t$ of the 
universal complex $(K^\bul_{\L}(\A),\D^\bul(x))$ of 
Theorem~\ref{thm:univcx}.  So $t \in \Sigma_{m}^{q}(M)$ if and only if 
$\dim H^{q}(K^\bul(\A),\D^\bul(t)) \ge m$.  An exercise in linear 
algebra shows that
\[
\Sigma_{m}^{q}(M)=\{t\in\T\mid
\rank \D^{q-1}(t) + \rank \D^{q}(t) \le \dim K^q(\A)-m\}.
\]

For $\la\in \C^{n}$, we have $\la\in \RR_{m}^{q}$ if $\dim 
H^{q}(A^{\bul},a_\la\wedge) \ge m$.  Denote the matrix of 
$a_\la\wedge: A^q(\A) \to A^{q+1}(\A)$ by $\mu^q(\la)$.  Then, as 
above, 
\[
\RR_{m}^{q}(A)=\{\la\in\C^{n}\mid 
\rank \mu^{q-1}(\la) + \rank \mu^{q}(\la) \le \dim A^{q}(\A)-m\}.
\]
Now $\dim A^{q}(\A)=\dim K^{q}(\A)=b_q(\A)$ and for each $q$, by 
Theorem~\ref{thm:approx}, we have 
$\rank\mu^{q}(\la)=\rank\D^{q}_{*}(\la)$.  Thus,
\begin{align*}
\Sigma_m^q(M)&=\{t\in\T\mid \rank \D^{q-1}(t) + \rank \D^{q}(t) \le 
b_q(\A)-m\},
\text{ and,}\\
\RR_{m}^{q}(A)&=\{\la\in\C^{n}\mid \rank \D^{q-1}_{*}(\la) + 
\rank \D^{q}_{*}(\la) \le b_q(\A)-m\},
\end{align*}
and the result follows.
\end{proof}

The cohomology support loci are known to be unions of 
torsion-translated subtori of $\T$, see \cite{Ar}.  In particular, all 
irreducible components of $\Sigma_{m}^{q}(M)$ passing through $\bone$ 
are subtori of $\T$.  Consequently, all irreducible components of the 
tangent cone are linear subspaces of $\C^{n}$.  So we have the 
following.

\begin{cor} \label{cor:falkconj}
For each $q$ and $m$, the resonance variety $\RR^{q}_{m}(A)$ is 
the union of an arrangement of subspaces in $\C^{n}$.
\end{cor}

\begin{rem} 
Several special cases of Theorem \ref{thm:tcone} and Corollary 
\ref{cor:falkconj} were previously known.  For $q=1$, these results 
were established by Cohen-Suciu~\cite{CS4}, see also 
Libgober-Yuzvinsky~\cite{Li,LY}.  For the discriminantal arrangements 
of Schechtman-Varchenko~\cite{SV}, they were established in \cite{C3}.  
In particular, as conjectured by Falk \cite[Conjecture~4.7]{Fa}, the 
resonance varieties $\RR_{m}^{q}(A)$ were known to be unions of linear 
subspaces in these instances.  Corollary \ref{cor:falkconj} above 
resolves this conjecture positively for all arrangements in all 
dimensions.
\end{rem}

\begin{rem} Theorem \ref{thm:tcone} and Corollary \ref{cor:falkconj} 
have been obtained recently by Libgober in a more general situation, 
see \cite{L2}.
\end{rem}

\section{Rational Local Systems} \label{sec:ratl}
Let $\A =\{H_{1},\dots,H_{n}\}$ be an arrangement of complex 
hyperplanes, and let $\la=(\la_1,\dots,\la_n)$ be a system of rational 
weights, $\la_j= k_j/N$, where $k_j\in\Z$ and $N\in\N$.  The 
representation $\rho:\pi_1(M)\to \C^*$, $\c_j\mapsto 
t_j=\exp(2\pi\ii\la_j)$, is unitary, and we call the associated local 
system $\LL$ on $M$ {\em rational}.  Since our primary interest is in 
this local system, we may assume without loss that the greatest common 
divisor of the integers $k_j$ is relatively prime to $N$.  In this 
section, we obtain combinatorial {\em upper} bounds on the local 
system Betti numbers $\dim_{\C} H^q(M;\LL)$.  We then use these bounds 
to study non-resonant local systems.

\subsection{Combinatorial Cohomology mod $N$} \label{subsec:OSmodN}
Let $A_\Q(\A)$ be the Orlik-Solomon algebra of $\A$, with rational 
coefficients, and generated by $\{a_1,\dots,a_n\}$.  If the underlying 
arrangement is clear, we write $A_\Q=A_\Q(\A)$.  Left-multiplication 
by the element $a_\la=\sum \la_j a_j \in A^1_\Q$ induces a 
differential on the Orlik-Solomon algebra, and we denote the resulting 
complex by $(A^{\bul}_\Q,a_\la\wedge)$.  Similarly, associated to the 
element $a_k=N a_\la=\sum k_j a_j$, we have the complex 
$(A^{\bul}_\Q,a_k\wedge)$.

\begin{lem} \label{lem:kequiv}
The complexes $(A^{\bul}_\Q,a_\la\wedge)$ and 
$(A^{\bul}_\Q,a_k\wedge)$ are chain equivalent.
\end{lem}
\begin{proof}
Define $\eta:A^{\bul}_\Q\to A^{\bul}_\Q$ by $\eta(a)=N^q a$ for $a\in 
A_\Q^q$.  Since $\eta$ is clearly an isomorphism, it is enough to show 
that $\eta$ is a chain map.  For $a$ as above, we have
\[
\eta(a_\la\wedge a)=N^{q+1}\sum \frac{k_j}{N}a_j\wedge a= N^q 
a_k\wedge a =
a_k\wedge\eta(a),
\]
and $\eta$ is a chain map.
\end{proof}

Denote the matrix of $a_\la\wedge:A^q_\Q \to A^{q+1}_\Q$ by 
$\mu^q(\la)$, and that of $a_k\wedge:A^q_\Q \to A^{q+1}_\Q$ by 
$\mu^q(k)$.  The entries of the latter are integers, so we consider 
the Orlik-Solomon algebra with integer coefficients, the associated 
complex $(A^{\bul}_\Z,a_k\wedge)$, and the reduction of this 
complex$\mod N$.  Let $(A^{\bul}_N,\bar{a}_k\wedge)$ be the reduction 
of $(A^{\bul}_\Z,a_k\wedge)\mod N$, where $A_N=A_{\ZN}$ denotes the 
Orlik-Solomon algebra with coefficients in the ring $\ZN$ and 
$\bar{a}_k=a_k\mod N$.  Denote the matrix of $\bar{a}_k\wedge:A^q_N 
\to A^{q+1}_N$ by $\bar\mu^q(k)$.

\subsection{Upper Bounds} \label{subsec:upper}
We now obtain combinatorial upper bounds on the the Betti numbers, 
$\dim_{\C} H^q(M;\LL)$, for a rational local system $\LL$.  We shall 
make use of the following elementary fact.

\begin{lem} \label{lem:zerodivisor}
Let $\zeta=\exp(2\pi\ii/N)$ be a primitive $N$-th root of unity, and 
let $f(z)\in\C[z^{\pm 1}]$ be a Laurent polynomial which satisfies 
$f(1)=0$, $f(\zeta)=0$, and $f'(1)\in\Z$.  Then $[f'(1)]$ is a 
zero-divisor in $\ZN$.
\end{lem}
\begin{proof}
Write $f(z)=z^{-m} p(z)$, where $p(z)\in\C[z]$ is a polynomial.  Since 
$p(\zeta)=f(\zeta)=0$, we have $p(z)=\Phi_N(z)\cdot q(z)$, where 
$\Phi_N(z)$ is the $N$-th cyclotomic polynomial and $q(z)\in\C[z]$.  
Since $p(1)=f(1)=0$ and $\Phi_N(1) \neq 0$, we also have $q(1)=0$.  
These considerations, and a brief calculation, reveal that 
$f'(1)=p'(1)=\Phi_N(1)\cdot q'(1)$.  Since $\Phi_N(1)$ divides $N$, we 
see that $[f'(1)]$ is a zero-divisor in $\ZN$.
\end{proof}

The main result of this section is the following.

\begin{thm}\label{thm:modN}
Let $\la=k/N$ be a system of rational weights, and let $\LL$ be the 
associated rational local system on the complement $M$ of $\A$.  Then, 
for each $q$,
\begin{equation} \label{eq:UB}
\dim_\C H^q(M;\LL) \le \rank_{\ZN} H^q(A^{\bul}_N,\bar{a}_k\wedge).
\end{equation}
\end{thm}
\begin{proof}
For any system of rational weights $\la=\frac{1}{N}(k_1,\dots,k_n)$, 
the rational local system $\LL$ is determined by the unitary 
representation $\rho:\pi_1(M)\to\C^*$ defined by $\c_j\mapsto 
\zeta^{k_j}$, where $\zeta=\exp(2\pi\ii/N)$.  This representation 
factors through the integers $\Z=\langle z \rangle$ as follows.  
Define $\xi:\pi_1(M) \to \Z$ by $\xi(\c_j)= z^{k_j}$ and define 
$\chi:\Z \to \C^*$ by $\chi(z)=\zeta$.  We then have $\rho = \chi 
\circ \xi$.

The cohomology of $M$ with coefficients in $\LL$ is isomorphic to that 
of the complex $(K^\bul(\A),\D^\bul(t))$, the specialization of the 
universal complex $(K^\bul_{\L}(\A),\D^\bul(x))$ at the point 
$t=(t_1,\dots,t_n)\in\T$, where $t_j=\exp(2\pi\ii k_j/N)=\zeta^{k_j}$.  
Thus for a rational local system, via the map $\L\to\C[z^{\pm 1}]$ 
defined by $t_j \mapsto z^{k_j}$, the specialization map $\L \to \C$ 
factors through the ring of Laurent polynomials in the single variable 
$z$.  To emphasize the dependence of the boundary maps of $K^\bul(\A)$ 
on $\zeta$, we shall write $\D^\bul(t)=\D^\bul(\zeta)$.  
By virtue of the above factorization, for each $q$, the matrix 
$\D^q(\zeta)$ may be realized as the specialization, $z \mapsto \zeta$,  
of a matrix $\D^q(z)$ with entries in $\C[z^{\pm 1}]$.

By Theorem~\ref{thm:approx}, the Aomoto complex $(A^{\bul}_R(\A),a_y 
\wedge)$ is chain equivalent to the linearization, 
$(K^\bul_{R}(\A),\D^\bul_*(y))$, of the universal complex 
$(K^\bul_{\L}(\A),\D^\bul(x))$.  
Both this chain equivalence and the construction of the Aomoto complex  
involve choosing a basis for the Orlik-Solomon algebra of $\A$. 
Making these choices in a consistent 
manner, we may assume without loss that the two complexes coincide, 
$(K^\bul_{R}(\A),\D^\bul_*(y))=(A^{\bul}_R(\A),a_y \wedge)$.  
Specializing, $y\mapsto k$, yields the complex 
$(A_\C,a_k\wedge)=(K^\bul(\A),\D^\bul_*(k))$.  Recall that we denote 
the matrix of $a_k\wedge:A^q \to A^{q+1}$ by $\mu^q(k)$.  Since for 
each $j$, we have $t_j=\chi \circ \xi(x_j)=\chi(z^{k_j})=\zeta^{k_j}$, 
these boundary maps may be realized as the derivatives at $z=1$ of the 
maps $\D^q(z)$ noted above: $\mu^q(k)=\D^q_*(k)=\frac{d}{dz} 
\D^q(z)\bigr\rvert_{z=1}$.

Recall from Lemma~\ref{lem:linear} that the boundary maps of the 
Aomoto complex consist of integral linear combinations of the 
indeterminates $y_j$.  Thus the entries of the matrices 
$\mu^q(k)=\D^q_*(k)$ above are integers, and, as in Section 
\ref{subsec:OSmodN} above, we may consider the complex 
$(A^{\bul}_\Z,a_k\wedge)$.  The reduction, 
$(A^{\bul}_N,\bar{a}_k\wedge)$, of this complex mod $N$ has boundary 
maps given by the matrices $\bar\mu^q(k)$, the reductions mod $N$ of 
$\mu^q(k)=\D^q_*(k)$.  By Lemma~\ref{lem:zerodivisor}, if a minor of 
$\bar\mu^q(k)$ is a unit in $\ZN$, then the corresponding minor of 
$\D^q(\zeta)$ is non-zero.  Thus, $\rank_\C \D^q(\zeta) \ge 
\rank_{\ZN} \bar\mu^q(k)$, and the result follows.
\end{proof}

There are rational local systems $\LL$ (resp.,~weight vectors $\la$) 
for which the inequality \eqref{eq:UB} is strict, as illustrated by 
the following example.

\begin{exm} \label{exm:ustrict}
Let $\A=\{H_1,\dots ,H_8\}$ be a realization of the MacLane ($8_{3}$) 
configuration, with defining polynomial 
\[
Q(\A)=xy(y-x)z(z-x-\omega^{2}y)(z+\omega y)(z-x)
(z+\omega^{2}x+\omega y),
\]
where $\omega$ is a primitive third root of unity and the ordering of 
the hyperplanes of $\A$ corresponds to that of the factors of $Q(\A)$.  
Consider the rational weights
\[
\la(u,v)=
\frac{u}{3}(1,0,2,1,2,2,1,0)+\frac{v}{3}(2,2,2,1,1,0,0,1),
\]
where $u,v\in\{0,1,2\}$, the associated points $t(u,v)=\exp(2\pi\ii 
\la(u,v))$ in the complex torus $(\C^{*})^{8}$, and the corresponding 
rational local systems $\LL=\LL(u,v)$ on the complement $M$ of $\A$.

As observed by Matei-Suciu \cite[Example~5.9]{MS}, the mod 3 
resonance variety, $\RR^{1}_{1}(A,\Z_3)$, of the Orlik-Solomon 
algebra of $\A$ contains a two-dimensional component.  In our 
notation, this component may be described by the equations
\[
C=\{\bar{k}=(k_{1},\dots,k_{8}) \in (\Z_3)^{8} \mid 
k_{1}+k_{5}=k_{2}+k_{8}= k_{3}+k_{4}=k_{6}+k_{7}=0\}.
\]
The (nine) points of $C$ may be realized as the reductions mod 3 of 
the integral weight vectors $k(u,v)=3\la(u,v)$ for $u,v\in\{0,1,2\}$.  
Thus for $k=k(u,v)$ with $(u,v)\neq (0,0)$, we have $\rank_{\Z_3} 
H^1(A_3,\bar{a}_k\wedge)=1$, as may be readily checked by direct 
calculation in the Orlik-Solomon algebra.

On the other hand, for each of the non-trivial rational local systems 
$\LL=\LL(u,v)$, we have $H^{1}(M;\LL)=0$.  Using the braided wiring 
diagram for this arrangement recorded in \cite[Example 8.6]{CS3}, one 
may compute the braid monodromy of $\A$, and the ensuing braid 
monodromy presentation of the fundamental group of $M$, see 
\cite{CS2}.  Then an exercise with this presentation and the Fox 
calculus reveals that the first local system cohomology of $M$ is 
trivial, $H^1(M;\LL)=0$, for each of the local systems $\LL=\LL(u,v)$.  
Thus, for each of the weight vectors $\la(u,v)$, the inequality 
\eqref{eq:UB} is strict.

Note that, since $H^1(M;\LL)=0$ for $\LL=\LL(u,v)$, we have 
$H^1(A^\bul_\C;a_{\la+m}\wedge)=0$ for every integral translate 
$m\in\Z^8$ of each of the weight vectors $\la=\la(u,v)$.
\end{exm}

\begin{exm} \label{exm:bothstrict}
Using the previous examples, one can construct arrangements and weight 
vectors for which both inequalities \eqref{eq:LB} and \eqref{eq:UB} 
are simultaneously strict.

Let $\A'$ be a generic section in $\C^2$ of the Ceva(3) arrangement 
from Example~\ref{exm:ceva}, and let $\A''$ be a generic section 
in $\C^2$ of the Maclane arrangement above.  Consider the weight 
vectors $\la'=\frac{1}{3}k'$ for $\A'$ and 
$\la''=\frac{1}{3}k''$ for $\A''$, where
\[
k'=(1,1,1,1,1,1,-2,-2,-2)\in\Z^9\quad\text{and}\quad
k''=(1,0,-1,1,-1,-1,1,0)\in\Z^8.
\]
Note that $\la''$ is an integral translate of the weight vector 
$\la(1,0)$ from the previous example and that $k''\equiv k(1,0)\mod 
3$.  Let $\LL'$ and $\LL''$ denote the local systems on $M(\A')$ and 
$M(\A'')$ corresponding to $\la'$ and $\la''$.  Since these 
arrangements are generic sections, our previous calculations and Euler 
characteristic arguments yield the Orlik-Solomon algebra cohomology 
and local system cohomology.

We record this information using Poincar\'e polynomials.  Write
\begin{align*}
P(A^\bul_\C(\A),\la,t)&=\sum \dim_\C H^i(A^\bul_\C(\A),a_\la\wedge) 
t^i,\\
P(M(\A),\LL,t)&=\sum \dim_\C H^i(M(\A),\LL) t^i,\ \text{and}\\
P(A^\bul_N(\A),k,t)&=\sum\rank_{\Z_N} 
H^i(A^\bul_N(\A),\bar{a}_k\wedge) t^i.
\end{align*}
Then, for $\A'$, we have
\[
P(A^\bul_\C(\A'),\la',t)=t+17t^2 \quad\text{and}\quad
P(M(\A'),\LL',t)=2t+18t^2.
\]
Furthermore, a computation in the Orlik-Solomon algebra reveals that 
\[
P(A^\bul_3(\A'),k',t)=P(M(\A'),\LL',t).  
\]
For $\A''$, we have
\[
P(A^\bul_\C(\A''),\la'',t)=P(M(\A''),\LL'',t)=13t^2 
\quad \text{and}\quad  P(A^\bul_3(\A''),k'',t)=t+14t^2.
\]

Now let $\A=\A' \times \A''$ be the product arrangement in $\C^4$, and 
let $k=(k',k'')\in\Z^{17}$.  Consider the weight vector 
$\la=\frac{1}{3}k=(\la',\la'')=\frac{1}{3}(k',k'') \in\C^{17}$, and 
the associated local system $\LL$ on $M(\A)$.  By construction, $\LL$ 
is the product local system on $M(\A)\cong M(\A') \times M(\A'')$.  
Similarly, the boundary map of the Orlik-Solomon algebra complex is 
compatible with the product structure $A(\A)\simeq A(\A') \times 
A(\A'')$.  Consequently, we may use the K\"unneth formula to obtain
\begin{align*}
P(A^\bul_\C(\A),\la,t)&=P(A^\bul_\C(\A'),\la'',t)\cdot 
P(A^\bul_\C(\A''),\la'',t)=13t^3+221t^4,\\
P(M(\A),\LL,t)&=P(M(\A'),\LL',t) \cdot 
P(M(\A''),\LL'',t)=26t^3+234t^4,\ \text{and}\\
P(A^\bul_3(\A),k,t)&=P(A^\bul_3(\A'),k',t) \cdot 
P(A^\bul_3(\A''),k'',t)= 2t^2+46t^3+252t^4.
\end{align*}
Thus, 
\[
\dim_\C H^i(A^\bul_\C(\A),a_\la\wedge) < \dim_\C H^i(M(\A),\LL) <
\rank_{\Z_3} H^i(A^\bul_3(\A),\bar{a}_k\wedge)
\] 
for $i=3,4$, and both inequalities \eqref{eq:LB} and \eqref{eq:UB} are 
strict.
\end{exm}

\subsection{Resonance} \label{subsec:res}
We use the results obtained above to study resonance phenomena.  Let 
$\A$ be an essential arrangement of $n$ hyperplanes in $\C^\ll$, with 
complement $M$.  We call a system of weights $\la\in\C^n$ {\em 
non-resonant} if the Betti numbers of $M$ with coefficients in the 
associated local system $\LL$ are minimal.  Since the boundary maps 
$\D^q(t)$ of the complex $K^\bul(\A)$ computing the local system 
cohomology generically take on their maximal ranks, the set of 
non-resonant weights is open and dense in $\C^n$.  We denote this set 
by $\UU=\UU(\A)$.  It may be described as the set of all $\la\in\C^n$ 
for which the sum, $\Sigma(\A)=\sum_{q=0}^\ll \dim H^q(M;\LL)$, of the 
local system Betti numbers is minimal.

Recall that an edge is a nonempty intersection of hyperplanes.  An 
edge is called {\em dense} if the subarrangement of hyperplanes 
containing it is irreducible: the hyperplanes cannot be partitioned 
into nonempty sets so that after a change of coordinates hyperplanes 
in different sets are in different coordinates.  This is a 
combinatorially determined condition, see \cite{STV}.  For each edge 
$X$, define $\la_X=\sum_{X \subseteq H_j}\la_j$.  Let 
$\A_\infty=\A\cup H_\infty$ denote the projective closure of $\A$, the 
union of $\A$ and the hyperplane at infinity in $\CP^\ll$, see 
\cite{OT2}.  Consider the following sets:
\begin{align*}
\VV=\VV(\A)&=\{\la\in\C^n \mid \la_X \notin \Z_{>0}
\ \text{for every dense edge $X$ of $\A_\infty$}\},\\
\WW=\WW(\A)&=\{\la\in\C^n \mid \la_X \notin \Z_{\ge 0}
\ \text{for every dense edge $X$ of $\A_\infty$}\}.  
\end{align*}

\begin{thm} \label{thm:nonres}
If $\la\in \WW$, then $H^q(M;\LL)=0$ for $q\neq\ll$ and $\dim
H^\ll(M;\LL)=|e(M)|$.
\end{thm}

Since $\dim K^q(\A)=b_q(\A)=\dim_\C H^q(M;\C)$, we see that this 
theorem minimizes $\Sigma(\A)$.  It follows that $\WW\subset \UU$, 
proving Theorem \ref{nonres} of the Introduction.  There is an 
important difference between these statements, however.  Theorem 
\ref{nonres} is an existence result which does not reveal the 
structure of $\UU$, while Theorem \ref{thm:nonres} identifies a large 
subset of $\UU$.  Here we state the two ingredients of the proof of 
Theorem \ref{thm:nonres} more precisely.  Choose a degree one 
polynomial $\alpha_j$ with kernel $H_j\in \A$ for $1 \leq j \leq n$.  
The Brieskorn algebra, $B(\A)$, is the graded $\C$-algebra generated 
by 1 and the holomorphic 1-forms $\omega_j=d \alpha_j/\alpha_j$.  Let 
$\omega_\la=\sum_{j=1}^n \la_j \omega_j$.  Then 
$(B^\bul(\A),\omega_\la\wedge)$ is a subcomplex of the twisted de Rham 
complex used to calculate $H^*(M;\LL)$.  The results establishing 
Theorem~\ref{thm:nonres} are:
\begin{enumerate}
\item If $\la\in \VV$, then $H^*(M;\LL)\simeq 
H^*(B^\bul(\A),\omega_\la\wedge)$.  This is a result of 
Esnault-Schechtman-Viehweg \cite{ESV}, refined by 
Schechtman-Terao-Varchenko \cite{STV}, and obtained by using Deligne's 
work \cite{De}.  The condition $\la\in \VV$ is a monodromy condition, 
imposed on a normal crossing divisor obtained by resolution of 
singularities.

\item The cohomology, $H^*(B^\bul(\A),\omega_\la\wedge)$, of the 
Brieskorn algebra is isomorphic to that of the Orlik-Solomon algebra, 
$H^*(A^\bul(\A),a_\la\wedge)$, \cite{OT2}.  The latter was studied by 
Yuzvinsky \cite{Yuz}, who showed that if $\la \in \WW$, then 
$H^q(A^\bul(\A),a_\la\wedge)=0$ for $q\neq\ll$, from which it follows 
that $\dim H^\ll(A^\bul(\A),a_\la\wedge)=|e(M)|$.
\end{enumerate}

We give a proof of Theorem~\ref{nonres} which uses only 
Theorem~\ref{thm:modN} and Yuzvinsky's result stated above, and thus 
avoids use of resolution of singularities.

\begin{proof}[Proof of Theorem~\ref{nonres}]
The alternating sum of the dimensions of the cochain groups $K^q(\A)$ 
is the Euler characteristic, so it suffices to prove that there exists 
a system of weights with cohomology groups as stated in the theorem.  
Yuzvinsky's argument may be applied with weights in an arbitrary 
field, and requires only that $\la_X$ be an invertible element for all 
dense edges.

Choose a prime $p>n$ and let $\la = k/p$, where $k_j=1$ for every $j$.  
Then for every (dense) edge $X$, $k_X=|X|\leq n < p$.  Thus $k_X$ is a 
unit in $\Z_p$, and $H^q(A^\bul_p(\A),\bar{a}_k\wedge)=0$ for 
$q\neq\ll$.  So by Theorem~\ref{thm:modN}, for the rational local 
system $\LL$ corresponding to $\la$, we have $H^q(M,\LL)=0$ for 
$q\neq\ll$, and we are done.
\end{proof}

\begin{rem}
There are in general many non-resonant weights which do not lie in the 
set $\WW$.  For instance, we have $\la(1,0)\in\UU\setminus\WW$ in 
Example \ref{exm:ustrict}.  Turning to resonant weights, the 
arrangement of three lines through the origin in $\C^2$ with weights 
$\la=(\la_1,\la_2,\la_3)\in\VV$ satisfying $\la_1+\la_2+\la_3=0$ shows 
that $\VV \not \subseteq \UU$.  Examples~\ref{exm:lstrict} and 
\ref{exm:ceva} show that there are resonant weights $\la$ for which 
$\la+m\not \in \VV$ for any $m\in \Z^n$.  The following summarizes our 
current understanding of resonant weights:
\begin{enumerate}
\item 
If $\la+m \in \VV$ for some $m\in \Z^n$, then $H^*(M;\LL)\simeq
H^*(A_\C^\bul,a_{\la+m}\wedge)$ and we have an effective algorithm.

\item
Otherwise, we have only the inequalities 
\[
\dim_\C H^q(A_\C^{\bul},a_\la\wedge)\leq \dim_\C H^q(M;\LL)
\leq \dim_\C H^q(M;\C)
\]
for arbitrary weights, and 
\[
\dim_\C H^q(A_\C^{\bul},a_\la\wedge)\leq \dim_\C H^q(M;\LL)\leq 
\rank_{\ZN} H^q(A^{\bul}_N,\bar{a}_k\wedge)
\]
for rational weights.
\end{enumerate}
\end{rem}

\begin{ack}
Portions of this work were carried out in the Spring of 1999, when 
the first author served as an Honorary Fellow in the Department of 
Mathematics at the University of Wisconsin-Madison.  He thanks the 
Department for its hospitality, and for providing an exciting and 
productive mathematical environment.  He also thanks the College of 
Arts {\&} Sciences at Louisiana State University for granting the 
Research Fellowship which made this visit possible.
\end{ack}

\bibliographystyle{amsalpha}

\vskip -3.0015pt

\end{document}